\documentclass[11pt,reqno]{amsart}
\usepackage{xy,epsfig,amsfonts,amssymb}

\catcode `\@=11

\def\numberbysection{\@addtoreset{equation}{section}
\renewcommand{\theequation}{\thesection.\arabic{equation}}}
\numberbysection
\def\subsubsection{\@startsection{subsubsection}{3}%
  \normalparindent{.5\linespacing\@plus.7\linespacing}{-.5em}%
  {\normalfont\bfseries}}
\def\cal{\mathcal}

\def\T{\mathcal T}

\def\S{{\Sigma}}

\DeclareMathOperator{\Ker}{Ker}
\DeclareMathOperator{\Ima}{Im}
\DeclareMathOperator{\napo}{\mathcal{N}AP}
\DeclareMathOperator{\Prim}{Prim}
\DeclareMathOperator{\Hom}{Hom}
\DeclareMathOperator{\Vertex}{Vert}
\DeclareMathOperator{\End}{\mathcal{E}nd}
\DeclareMathOperator{\Free}{\mathcal{F}ree}

\def\nap{non-\-asso\-cia\-ti\-ve per\-mu\-ta\-ti\-ve}
\def\PL{\mathcal{PL}}
\def\RT{\mathcal{RT}}
\def\HO{\mathcal{HO}}

\xyoption{all}
\CompileMatrices

\begin{document}

\title{A rigidity theorem for pre-Lie algebras}
\author[M. Livernet]{Muriel Livernet}

\maketitle

\begin{center}
Universit\'e Paris 13, Institut Galil\'ee, LAGA \\
{\small Avenue Jean-Baptiste Cl\'ement, 93430
Villetaneuse, France \\
e-mail: livernet@math.univ-paris13.fr}
\end{center}

\date{\today}


\begin{abstract}

In this paper we prove a ``Leray theorem'' for pre-Lie algebras. We
define a notion of ''Hopf'' pre-Lie algebra: 
it is a pre-Lie algebra together with a nonassociative permutative
coproduct $\Delta$ and a compatibility relation 
between the pre-Lie product and the coproduct $\Delta$. A
nonassociative permutative algebra is a vector space together with a
product satisfying 
the relation $(ab)c=(ac)b$. 
A nonassociative permutative coalgebra is the dual notion. We prove
that any connected "Hopf" pre-Lie algebra is a free pre-Lie
algebra. It uses the description of pre-Lie algebras in term of rooted
trees developped by 
Chapoton and the author. We interpret also this theorem by way of cogroups in the category of pre-Lie algebras.\\

\end{abstract}


\section*{Introduction}

The classical Leray theorem (\cite{Le45}, \cite{MM65}) asserts that
any graded connected commutative and cocommutative Hopf algebra is a
free commutative algebra and free cocommutative coalgebra. Loday and
Ronco \cite{LR04} obtain a similar result for unital infinitesimal bialgebras: any
such object (graded connected) is a free associative algebra as well
as a free coassociative coalgebra. The same kind of result is also
obtained by Foissy in \cite{F05} for dendriform, codendriform
bialgebras. These three results are similar in the operadic framework:
given the operad $\mathcal P$, where $\mathcal P$ is the commutative operad
or associative operad or dendriform operad, then any graded connected $\mathcal
P$-algebra which is also a $\mathcal P$-coalgebra equipped with a relation between
the coalgebra and algebra structures--also called ``distributive
law''--is rigid in the sense that it is free as a 
$\mathcal P$-algebra and as a $\mathcal P$-coalgebra. In
this paper we obtain a result of that kind for pre-Lie algebras,
although the originality of our result is that the co-structure
involved is not pre-Lie, but \nap. In fact, such a result can be also
interpreted as a ``Cartier-Milnor-Moore'' theorem where the primitive
part is a vector space without any additional structure and the
envelopping algebra is the free $\mathcal P$-algebra functor. In that
sense, such a result exists for Zinbiel algebras as a consequence of a theorem by Ronco  \cite{R00} concerning
the primitive elements of a free dendriform algebra.

\medskip

PreLie algebras have been of interest since the works of Vinberg
\cite{V63} and Gerstenhaber \cite{G63}. 
In \cite{ChaLi01} using operad theory,  pre-Lie algebras are described in terms of rooted trees, linking this structure
to renormalisation theory \`a la Connes and Kreimer \cite{ConK98}. In
this paper we describe another structure based on rooted trees, called
\nap\  algebras (see e.g. \cite{DL02}). A {\sl\nap\  algebra} is a vector space together with a
bilinear product satisfying the relation $(ab)c=(ac)b$.
We state the following rigidity theorem

\medskip

\noindent {\bf Theorem--} \it Any  pre-Lie algebra, together with a \nap\ connected coproduct satisfying the distributive law 
$$\Delta(a\circ b)=\Delta(a)\circ b+a\otimes b$$ 
is a free pre-Lie algebra and a free \nap\  coalgebra.\rm

\medskip

The proof of the theorem involves an idempotent in the vector space of endomorphism of $L$, which
annihilates $L^2$, the space of decomposable elements (see the fundamental lemma \ref{fondlemma}).

Although this paper can be read without referring to operads, the ideas coming
from operad theory
are always present, especially in the last section where we link our
result to the theory of cogroups in the category of algebras over an
operad as
developped by
Fresse in \cite{Frbull} and  \cite{Fr98}.

\medskip

The paper is organized as follows: the first section is devoted to
material concerning operads, 
rooted trees and pre-Lie algebras; in the second section we introduce \nap\  algebras and coalgebras; in the
third section we state and prove the main theorem, assuming the fundamental lemma; the latter
is proved in the fourth section; finally, the fifth section concerns cogroups, 
where we study the special case of associative algebras and pre-Lie algebras.

\medskip

\noindent{\bf Notation.} The ground field $K$ is of characteristic $0$. The symmetric group on
$n$ elements is denoted by $\Sigma_n$. For any vector space $V$, $\Sigma_n$ acts
on $V^{\otimes n}$, on the left by
$$\sigma\cdot (v_1\otimes\ldots\otimes v_n)=v_{\sigma^{-1}(1)}\otimes\ldots\otimes v_{\sigma^{-1}(n)}.$$
For any subgroup $G$ of $\Sigma_n$ the subspace of $V^{\otimes n}$ of invariants under $G$
 to is denoted $(V^{\otimes n})^G=\{x\in V^{\otimes n}| \sigma\cdot x=x,\forall \sigma\in G\}.$

\section{pre-Lie algebras and rooted trees}

\subsection{Operads}\label{operads} In this section, we review the
material needed for this article. 
For a complement on algebraic operads we refer to Ginzburg and
Kapranov \cite {GK94} or Loday \cite {Lo96}.   

\smallskip

A {\sl $\Sigma$-module} $\mathcal M=\{\mathcal M(n)\}_{n>0}$
is a collection of (right) $\Sigma_n$-modules. Any $\Sigma_n$-module $M$
gives rise to a $\Sigma$-module $\mathcal M$ by setting $\mathcal
M(q)=0$ if $q\not =n$ and $\mathcal M(n)=M$.

\bigskip

An {\it operad} is a right $\Sigma$-module $\{\mathcal
O(n)\}_{n>0}$ such that $\mathcal O(1)=K$, together with 
composition products: for $1\leq i\leq n$

\smallskip

\begin{center}
\begin{tabular}{cccccc}
$\circ_i$:&$\mathcal O(n)$&$\otimes$&$\mathcal O(m)$&$\longrightarrow $&$\mathcal O(n+m-1)$ \\
&$p$&$\otimes$&$q$ &$\mapsto$& $p\circ_i q$
\end{tabular}
\end{center}

\smallskip

These compositions are subject to associativity conditions,
unitary conditions and equivariance conditions with respect to the
action of the symmetric group. 

\medskip

Any vector space $A$ yields an operad $\End_A$ with $\End_A(n)=\Hom(A^{\otimes n},A)$
where 
$$(f\circ_i g)(a_1,\ldots,a_{n+m-1})=f(a_1,\ldots,a_{i-1},g(a_i,\ldots,a_{i+m-1}),a_{i+m},\ldots,a_{n+m-1}).$$ 
An {\sl algebra over an operad ${\mathcal O}$}, or {\sl ${\mathcal O}$-algebra}, 
is a vector space $A$ together with an operad morphism from ${\mathcal O}$ to $\End_A$. This is equivalent
to some data $$ev_n: {\mathcal O}(n)\otimes_{\S_n} A^{\otimes n}\rightarrow A$$
satisfying associativity conditions with respect to the compositions.

One of the most important point in the theory of algebraic operads is the following:

\medskip


{\it The free ${\mathcal O}$-algebra generated by $V$ is the space
\begin{equation}\label{freeOalgebra}
{\mathcal O}(V)=\oplus_{n\geq 1} {\mathcal O}(n)\otimes_{\S_n} V^{\otimes n}
\end{equation}
where the maps $ev_n$ are described in terms of the $\circ_i's$.}

\medskip

The following objects are defined as usual: given a $\S$-module $M$,
the {\sl free operad generated by $M$}, 
denoted by $\Free(M)$, exists and satisfy the usual universal
property. 
The notion of an {\sl ideal} of an operad exists: when modding out an
operad by an ideal, one gets an operad. 
The operads we are concerned with are {\sl quadratic binary operads}
that is they are of the form
$\Free(M)/<R>$ where $M$ is a $\S_2$-module, 
$R$ is a sub-$\S_3$-module of $\Free(M)(3)$ and $<R>$ is the ideal
of $\Free(M)$ generated by $R$. The space $M=\Free(M)(2)$ is the space
of ``operations" whereas the space $R$ is the space of ``relations". For
instance the operad defining associative algebras 
is $\Free(K[\S_2])/<R_{As}>$ where $K[\S_2]$ is the regular
representation of $\S_2$ and $R_{As}$ is the free $\S_3$-module generated by 
$\mu\circ_1\mu-\mu\circ_2\mu$, where $\mu$ is a generator of $K[\S_2]$.


\subsection{Rooted trees}A {\sl rooted tree} is a nonempty connected graph without loop
together with a distinguished vertex called {\sl the root}. This root gives an orientation of the graph:
edges are oriented towards the root. The set of vertices of a tree $T$
is denoted by $\Vertex(T)$. The orientation induces a partial order
on $\Vertex(T)$, the root being the minimal element. The {\sl degree} of a tree is the number of vertices.
A {\sl $n$- labeled rooted tree} is a rooted tree together with a
bijection between $\Vertex(T)$ and the set  $\{1,\ldots,n\}$. A
{\sl $n$-heap-ordered tree} is a $n$-labeled rooted tree where the
bijection respects the partial order on $\Vertex(T)$. 
The space $\RT(n)$ is the vector space spanned by the $n$-labeled rooted
trees, and the space $\HO(n)$ is the one spanned by the $n$-heap-ordered
trees. The right action of $\S_n$ on $\RT(n)$ is the action on the
labeling.
For a vector space $V$ denote by $\RT(V)$ the space
$\sum_n\RT(n)\otimes_{\S_n} V^{\otimes n}$. If $V$ is equipped 
with a basis $\{e_\alpha\}_\alpha$ then the space $\RT(V)$ is
the vector space spanned by the rooted trees labeled by the set
$\{e_\alpha\}_\alpha$.

Following the notation
of Connes and Kreimer \cite{ConK98}, any tree $T$ writes 
$$T:=B(r,T_1,\ldots,T_k)$$
where $r$ is the root (or the labeling of the root) and $T_1,\ldots, T_k$ are trees. The {\sl arity of $T$} is $k$, the number of incoming edges of its root. Note that the list of the trees $T_1,\ldots, T_k$ is unordered.


\subsection{PreLie algebras} A {\sl (right) pre-Lie}-algebra is a vector space $L$ together with a
product $\circ$ satisfying the relation
\begin{equation}\label{relprelie} 
(x\circ y)\circ z-x\circ(y\circ z)=(x\circ z)\circ y-x\circ(z\circ y)
\end{equation}
These algebras appeared also under the name  {\sl right-symmetric algebras}
or {\sl Vinberg algebras} (if we deal with the left relation). The
terminology of pre-Lie algebra is the one of Gerstenhaber \cite{G63}. A pre-Lie algebra $L$ yields
a Lie algebra structure on $L$ with the bracket $[a,b]=a\circ b-b\circ a$, which is denoted by $L_{Lie}$.

\medskip

A theorem by Chapoton and the author links pre-Lie
algebras to rooted trees as follows:


\subsection{Theorem}\cite{ChaLi01}\label{chaliv}\it The vector space $\RT(V)$ endowed with the product
$$S\circ T=\sum_{v\in\Vertex(S)} S\circ_v T,$$
where $S\circ_v T$ is the tree obtained by grafting the root of $T$ on the vertex $v$ of $S$,
is the free pre-Lie algebra generated by $V$. \rm

\medskip

The proof of this fact uses the description of pre-Lie algebras in terms
of algebras over an operad as well as the relation (\ref{freeOalgebra}). One can find some different proofs in Dzhumadil{$'$}daev and L{\"o}fwall \cite{DL02} or in Guin and Oudom \cite{GO04}.


\subsection{Remark}\label{trick} Using the pre-Lie product in $\RT(V)$, one has:\begin{multline}\label{trickformula}
B(v,T_1,\ldots,
T_{n-1},T_n)=B(v,T_1,\ldots,
T_{n-1})\circ T_n\\
-\sum_{i=1}^{n-1}B(v,T_1,\ldots,T_i\circ T_n,\ldots,
T_{n-1})
\end{multline}
Hence a proposition concerned with the pre-Lie product in $\RT(V)$ can be proved 
by a double induction, on the degree and on the arity of a tree $T$.


\section{Nonassociative permutative algebras and coalgebras}


\subsection{Definition}

\label{napdef} 
A  {\sl \nap\  algebra} is a vector space $L$ equipped with a bilinear product satisfying the relation
$$(ab)c=(ac)b,\ \forall a,b,c \in L.$$
Note that it is indeed the definition of a ``right" \nap\  algebra and we can define
what is a ``left" \nap\  algebra. Note also that a {\sl Novikov} algebra 
is precisely
a (left) pre-Lie algebra whose product satisfies the (right) \nap\ 
relation, see \cite{Os92} for instance. 
An associative algebra which is a \nap\  algebra is a
permutative algebra in the terminology of Chapoton in \cite{Cha01}.


\subsection{Definition} Let $\cal T$ be the following operad: $\cal T(n)$ is the $\S_n$-module
of $n$-rooted
trees $\RT(n)$; for trees $T\in\RT(n)$ and $S\in \RT(m)$, 
the composition $T\circ_i S$ is the rooted tree obtained by substituting the tree $S$ for
the vertex $i$ in $T$~: the outgoing edge of $i$, if exists, becomes the
outgoing edge of the root of $S$; incoming edges of $i$ are grafted on
the root  of $S$. Then, it is easy to check that 
these compositions endow $\cal T$ with a structure of an operad. Example:

$$\vcenter
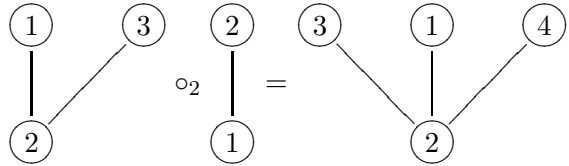
{\xymatrix{*++[o][F-]{1} \ar@{-}[d] & 
*++[o][F-]{3} \ar@{-}[dl]\\
*++[o][F-]{2}}} \circ_2
\vcenter{
\xymatrix{*++[o][F-]{2} \ar@{-}[d] \\
*++[o][F-]{1}}}=\vcenter{\xymatrix{*++[o][F-]{3}\ar@{-}[dr]&*++[o][F-]{1}\ar@{-}[d]&*++[o][F-]{4}\ar@{-}[dl]\\
&*++[o][F-]{2}}}$$

\medskip

If we compare with the operad $\PL$ based on rooted trees defining pre-Lie algebras in \cite{ChaLi01}, this 
composition is a summand of the composition in $\PL$.

Note that the operad $\T$ can also be defined on the set of $n$-labeled rooted trees: it is
in fact an operad in the category of sets.


\subsection{Proposition}\label{opnap}\it The previous operad $\cal T$ is the operad defining \nap\  algebras. \rm

\medskip

\noindent {\sl Proof--} Following the notation of section \ref{operads}, the operad defining \nap\  algebras is the operad $\napo=\Free(K[\S_2])/
<R_{nap}>$ where $R_{nap}$ is the sub-$\S_3$-module of $\Free(K[\S_2])(3)$
generated by the element $(\mu\circ_1\mu)\cdot (Id-\tau_{23})$ with
$\tau_{23}$ being the transposition $(132)$. There is a morphism of operads from
$\napo$ to $\cal T$ sending $\mu$ to the tree 
$\def\objectstyle{\scriptstyle}
\def\labelstyle{\scriptstyle}
\vcenter{\xymatrix@-1.5pc{
*++[o][F-]{2}\ar@{-}[d]\\
*++[o][F-]{1}}}.$ 
This morphism is well defined since the tree
$\def\objectstyle{\scriptstyle}
\def\labelstyle{\scriptstyle}
\vcenter{\xymatrix@-1.5pc{
*++[o][F-]{2}\ar@{-}[d]\\
*++[o][F-]{1}}}\circ_1 \vcenter{\xymatrix@-1.5pc{
*++[o][F-]{2}\ar@{-}[d]\\
*++[o][F-]{1}}}=\vcenter{\xymatrix@-1.5pc{
*++[o][F-]{2}\ar@{-}[d]& *++[o][F-]{3} \ar@{-}[dl]\\
*++[o][F-]{1}}} $
is invariant under the permutation $\tau_{23}$.

\smallskip

The inverse morphism is defined by induction on the number of vertices of $T$.
Any tree can be written $T=B(i,T_1,\ldots,T_n)$ where $i$ is the root of $T$, thus
$T=B(i,T_1)\circ_i B(i,T_2,\ldots, T_n)$ (up to an appropriate permutation). It is easy to check
that this decomposition does not depend on the choice of $T_1$. Hence by induction, 
it follows that we define a morphism of operads inverse to the previous one.\hfill$\Box$


\subsection{Corollary}\label{napfree} \it Let $V$ be a vector
space. The free \nap\  algebra
generated by $V$ is the vector space of rooted trees $\RT(V)$, 
endowed with the following product:
let $T=B(i,T_1,\ldots, T_n)$ and $S$ be two trees in $\RT(V)$ 
then $T\cdot S=B(i,T_1,\ldots,T_n,S)$. \rm

\medskip

\noindent{\sl Proof--} From relation (\ref{freeOalgebra}) and proposition \ref{opnap} the free \nap\  algebra on $V$ is 
$$\napo(V)=\bigoplus_{n\geq 1} \RT(n)\otimes_{\S_n} V^{\otimes n}$$
where the product is given by the composition of the operad, i.e. 

$$T\cdot S=\def\objectstyle{\scriptstyle}
\def\labelstyle{\scriptstyle}
(\vcenter{\xymatrix@-1.5pc{
*++[o][F-]{2}\ar@{-}[d]\\
*++[o][F-]{1}}}\circ_2 S)\circ_1 T=\vcenter{\xymatrix@-1.5pc{
*++[o][F-]{S}\ar@{-}[d]\\
*++[o][F-]{1}}}\circ_1 T,$$
which is the tree obtained by grafting $S$ on the root of $T$.\hfill $\Box$

\subsection{Remark} As pointed out before, since the operad $\T$ can be defined in the category of
sets, the corollary is also true when replacing \nap\  algebras with ``right-commutative magma"  in
the terminology of Dzhumadil{$'$}daev and L{\"o}fwall \cite{DL02}: these two results
were proved in their paper, using different methods.


\subsection{Definitions}\label{defconnexe} Let $(C,\Delta)$ be a
vector space $C$ together with a coproduct
$\Delta: C\rightarrow C\otimes C$. 
The following defines a filtration on $C$:

\medskip

$$\begin{tabular}{rl}
$\Prim C=C_1=$&$ \{x \in C | \Delta(x)=0\}$ \\
$C_n=$& $\{x\in C| \Delta(x)\in \sum_{i=1}^{n-1}
C_i\otimes C_{n-i}\}$
\end{tabular}$$

\medskip

The vector space $(C,\Delta)$ is said to be {\sl connected} if 
$C=\bigcup_{n\geq 1} C_n$. Note that  any graded $(C,\Delta)$  such that $C_0=0$ is connected.

Define $\Delta^k: C\rightarrow C^{\otimes (k+1)}$ by

$$\begin{tabular}{rl}
$\Delta^0=$& $Id$\\
$\Delta^1=$&$\Delta$ \\
$\Delta^{k+1}=$& $(\Delta\otimes Id^{\otimes k})\Delta^k=(\Delta^k\otimes Id)\Delta$
\end{tabular}$$

We use Sweedler notation for $\Delta^k$:  

$$\Delta^k(x)=\sum x_{(1)}\otimes\ldots\otimes x_{(k+1)},\ \forall x\in C.$$

\subsection{Lemma}\label{connected}\it Let $(C,\Delta)$ be a vector
space together with a coproduct; let $P_n$ be the vector space of cooperations from $C$ to $C^{\otimes (n+1)}$
built on $\Delta$. More explicitly: $P_0=span\{Id\}$; $P_1=span\{\Delta\}$;  $O\in P_{n}$ if and only if there exists  
$0\leq m\leq n-1$ and $O_1\in P_{m},\ O_2\in P_{n-m-1}$ such that 
$O=(O_1\otimes O_2)\Delta$.

\smallskip

\begin{center}
If $x\in C_n$ then $O(x)=0,\ \forall O\in P_n$. 
\end{center}
\rm

\medskip

\noindent{\sl Proof--} The lemma is trivial for $n=1$. Let $x$ be in $C_n$ and $O=(O_1\otimes O_2)\Delta$ in 
$P_n$, with $O_1\in P_m$. Since $\Delta(x)=\sum x_i\otimes y_{n-i}$,
with $x_i\in C_i$ and $y_{n-i}\in C_{n-i}$, we get by induction: $O_1(x_i)=0$ if $m\geq i$ and $O_2(y_{n-i})=0$
if $n-m-1\geq n-i$, i.e. if $m\leq i-1$; thus $(O_1\otimes O_2)(\Delta(x))=0$.\hfill$\Box$


\subsection{Definition} A {\sl \nap\  coalgebra} is a vector space together
with a coproduct $\Delta: C\rightarrow C\otimes C$ satisfying 
$(1-\tau_{23})(\Delta\otimes id)\Delta=0$. This is the dual notion of
a \nap\  algebra.

\medskip

The following lemma is an immediate consequence of the definition of \nap\  coalgebras.

\subsection{Lemma}\label{deltaklemma}\it 
If $C$ is a \nap\  coalgebra then 
$$\Ima(\Delta^k)\subset (C^{\otimes (k+1)})^{\S_1\times\S_k}.$$\rm


\subsection{Theorem}\label{freeconap} \it Let $V$ be a vector space. Then $\RT(V)$ together with the coproduct

$$\Delta(B(v,T_1,\ldots,T_n))=\sum_{i=1}^n B(v,T_1,\ldots,\widehat{T_i},\ldots,T_n)\otimes T_i$$
is the free \nap\  connected coalgebra. \rm

\medskip

\noindent {\sl Proof--} It is the dual statement of corollary \ref{napfree}, since $\RT(V)$ with the coproduct
is the graded dual of $\RT(V)$ considered as the free \nap\  algebra on $V$. It is graded by the number of vertices, hence
connected.  \hfill$\Box$


\section{The space of rooted trees and the main theorem}


\subsection{Definition}\label{defaction} Let $L$ be a pre-Lie algebra and $M$ be a vector
space. $M$ is a {\sl right $L$-module} if there exists a map $\circ:  M\otimes
L\rightarrow M$ such that
$$(m\circ l_1)\circ l_2-(m\circ l_2)\circ l_1=m\circ [l_1,l_2].$$
Equivalently: $M$ is a right $L$-module if and only if $M$ is a $L_{Lie}$-module in the terminology of Lie algebras.

Examples: $L$ is a $L$-module; $L^{\otimes n}$ is a $L$-module via
derivation as pointed out by Guin and Oudom in \cite{GO04}:
$$(x_1\otimes\ldots\otimes x_n)\circ y=\sum_{i=1}^n x_1\otimes\ldots\otimes
x_i\circ y\otimes\ldots\otimes x_n.$$


\subsection{Proposition}\it Consider the vector space $\RT(V)$ 
endowed with its pre-Lie product $\circ$ and its \nap\  coalgebra product $\Delta$. The
following relation is satisfied

\begin{equation}\label{dlaw}
\begin{split}
\Delta(x\circ y)=&x\otimes y+ \Delta(x)\circ y \\
=&x\otimes y+x_{(1)}\circ y\otimes x_{(2)}+x_{(1)}\otimes x_{(2)}\circ
y\\
\end{split}
\end{equation}\rm

\medskip

\noindent{\sl Proof--} For $S=B(v,S_1,\ldots,S_n)$ and $T$ in
$\RT(V)$, one has
\begin{equation*}
\begin{split}
\Delta(S\circ T)=&\Delta(B(v,S_1,\ldots,S_n,T))+\sum_{i=1}^n
\Delta(B(v,S_1,\ldots,S_i\circ T,\ldots, S_n))\\
=& S\otimes T+\sum_{i=1}^n
B(v,S_1,\ldots,\widehat{S_i},\ldots,S_n,T)\otimes S_i\\
&+\sum_{i\not=j} B(v,\ldots,\widehat{S_i},\ldots,S_j\circ
T,\ldots,S_n)\otimes S_i\\
&+\sum_i
B(v,S_1,\ldots,\widehat{S_i},\ldots,S_n)\otimes S_i\circ T\\
=&S\otimes T+\Delta(S)\circ T.\hspace{7cm} \Box
\end{split}
\end{equation*}


\subsection{Corollary}\label{propunicitemor}\it Let $H$ be a pre-Lie algebra, $V$ a vector
space and $\phi:\PL(V)\rightarrow H$ a morphism of pre-Lie algebras.
There exists a unique application $L: \PL(V)\rightarrow H\otimes H$
such that 

\begin{equation}\label{relunicitemor}
\begin{tabular}{cll}
$L(v)$&=\ \ 0&\ $\forall v\in V,$\\
$L(a\circ b)$&$=\phi(a)\otimes\phi(b)+L(a)\circ \phi(b)$&\ $\forall a,b\in\PL(V).$
\end{tabular}
\end{equation}
\rm

\medskip

\noindent{\sl Proof--} The existence is given by
$L=(\phi\otimes\phi)\Delta$ and the relation (\ref{dlaw}).
To prove unicity it is sufficient to prove the following: if
$L:\PL(V)\rightarrow H\otimes H$ satisfies $L(a\circ
b)=L(a)\circ\phi(b)$ and $L(v)=0$ then $L=0$. This is proved by
induction on the number of vertices of a tree $T$ and the arity of
$T$, thanks to the  relation (\ref{trickformula}). \hfill$\Box$


\subsection{Main Theorem}\label{maintheorem} \it Let $(H,\circ_H,\Delta_H)$ be a vector space
together with a pre-Lie product and a \nap\  connected coproduct satisfying the
relation $\Delta_H(a\circ_H b)=a\otimes b+\Delta_H(a)\circ_H b$. There is an
isomorphism
of pre-Lie algebras and
of \nap\  coalgebras between
$H$ and $(\RT(\Prim H),\circ,\Delta)$.\rm

\medskip

The proof of this theorem is similar to the proof of the rigidity of
unital infinitesimal bialgebras of Loday-Ronco \cite{LR04} and
the one of dendriform Hopf algebras of Foissy \cite{F05}. It relies
mainly on the fundamental lemma \ref{fondlemma}. To state this lemma
we need the following definition


\subsection{Definition--Notation}\label{lesak} Let $(H,\mu)$ be a
pre-Lie algebra. We define linear operators $A_k: H^{\otimes
  k}\rightarrow H$ by induction on $k$:

$$\begin{tabular}{rl}
$A_1=$&$Id$, \\
$A_2=$&$\mu$, \\
$A_{k+1}=$&$\sum_{l=1}^{k} {k-1\choose l-1} \mu(A_l\otimes A_{k+1-l}).$
\end{tabular}$$


\subsection{Remark} The operators $A_k$ can be viewed as elements of
$\RT(k)$ as follows. Define the scalar product on the basis of rooted
trees by 
$$<S,T>=\begin{cases} 1 & {\rm if\ } S=T \\ 0& {\rm if\  not.\
  }\end{cases}$$ Let $\HO(k)$ be the space of $k$-heap-ordered trees. Then
 
$$A_k=\sum_{U\in\HO(k)} c(U) U$$
where $c(U)$ is defined by induction:
$$c(U)=\sum_{l=1}^k {k-1\choose l-1} \sum_{T\in\HO(l),T'\in\HO(k+1-l)}
  c(T)c(T') <T\circ T',U>.$$


\subsection{Fundamental Lemma}\label{fondlemma}\it  Let
$(H,\circ,\Delta)$ be a pre-Lie algebra, \nap\  connected coalgebra satisfying
the relation (\ref{dlaw}). The linear morphism

\begin{equation}\label{defidempotent}
\begin{tabular}{cccl}
$e:$& $H$& $\rightarrow$ & $H$ \\
& $x$ & $\mapsto$ & $x+\sum_{k\geq 1} \frac{(-1)^k}{k!} A_{k+1}\Delta^k(x)$
\end{tabular}
\end{equation}

satisfies the following properties: 

\begin{itemize}
\item[1)] $e$ is a projector onto $\Prim(H)$.
\item[2)] $e(x\circ y)=0$,\ $\forall x,y \in H.$
\end{itemize}
\rm


\subsection{Corollary}\label{cor1} \it Let $(H,\mu)$ be a pre-Lie algebra together
with a \nap\  connected coproduct satisfying relation (\ref{dlaw}). Then
$$H=\Prim H\oplus H^2,$$ 
where $H^2=\Ima(\mu)$.\rm

\medskip

\noindent{\sl Proof--} By lemma \ref{fondlemma}: property (1) implies
that forall  $x\in H$ one has
$x=e(x)+x-e(x)$, where $e(x)\in\Prim H$ and $x-e(x)=-\sum_{k>0} \frac{(-1)^k}{k!} A_{k+1}\Delta^k(x)=\mu(y,z)$.
Property (2) implies that $\Prim H\cap H^2=0$.\hfill$\Box$


\subsection{Corollary}\label{cor2}\it  Let $(H,\mu)$ be a pre-Lie algebra together
with a \nap\  connected coproduct satisfying relation (\ref{dlaw}). Then $H$ is generated
as a pre-Lie algebra by $\Prim(H)$.\rm

\medskip

\noindent{\sl Proof--} Following definition \ref{defconnexe}, since
$H$ is connected, it decomposes as
$H=\bigcup_{n\geq 1}
H_n$; any $x\in H_n$ satisfies
$\Delta^k(x)\in \sum_{(i_1,\ldots,i_{k+1})} H_{i_1}\otimes\ldots\otimes H_{i_{k+1}}$
where $i_j<n,\forall j$. The proof is performed by induction on $n$:
for $n=1$, it is true since $H_1=\Prim H$.
If $x\in H_n$ then $x=e(x)-\sum_{k>0} \frac{(-1)^k}{k!} A_{k+1}\Delta^k(x)$.
But $e(x)\in \Prim H$ and by induction, each component in $\Delta^k(x)$ is generated by elements
of $\Prim(H)$, so is $A_{k+1}\Delta^k(x)$.\hfill $\Box$


\subsection{Proof of the main theorem} Set $V:=\Prim(H)$. 
Denote by $i:V\rightarrow H$ and $\iota: V\rightarrow \RT(V)$
the natural injections.

\smallskip

Since $(\RT(V),\circ)$ is the free
pre-Lie algebra on $V$, there is a morphism of pre-Lie algebras $\phi:\RT(V)\rightarrow H$, such that
the diagram
$$\xymatrix{ V  \ar[r]^{i} \ar[d]_{\iota} & H \\
{\mathcal RT}(V) \ar@{-->}[ur]_{\phi}& }$$
commutes. The pre-Lie algebra $H$ is generated by $V$ by  corollary \ref{cor2}, then $\phi$
is surjective.

\smallskip

Consider the two linear maps $\Delta_H\phi$ and $(\phi\otimes\phi)\Delta$ from  $\RT(V)$
to $H\otimes H$. An easy computation proves that these two maps satisfy the conditions (\ref{relunicitemor})
of  corollary \ref{propunicitemor}. Then by unicity, these two maps coincide, thus
$\phi$ is a morphism of \nap\  coalgebras.

\smallskip

Assume $\Ker\phi\not=0$ and let $a\in\Ker\phi$ of minimal degree
$n$. Since $(\phi\otimes\phi)\Delta(a)=\Delta_H\phi(a)=0$ the element
$a$ lies in $\Ker(\phi)\otimes\RT(V)+\RT(V)\otimes\Ker(\phi)$. Since $a$
is of minimal degree it implies that $a=0$, thus a
contradiction. Hence $\phi$ is
injective. \hfill $\Box$

\section{Proof of the fundamental lemma \ref{fondlemma}}


Assume $(H,\mu,\Delta)$ satisfies the hypothesis of lemma
\ref{fondlemma}. We define linear maps $U_{k+1}: H^{\otimes (k+1)}\rightarrow H^{\otimes 2}$, for $k>0$ by $U_{k+1}=\sum_{l=1}^{k}
  {k-1 \choose l-1} A_l\otimes A_{k+1-l}$, where the $A_i's$ were defined in \ref{lesak}. 
  Hence $A_k=\mu U_k$ for $k>1$.


\subsection{Proposition}\it Assume $k\geq 1$. On the vector space 
$$Z_{k+1}= (H^{\otimes (k+1)})^{\S_1\times
\S_k}\cap\{x\in H^{\otimes (k+1)} | (\Delta\otimes Id^{\otimes k})(x)\in (H^{\otimes (k+2)})^{\S_1\times
\S_{k+1}}\}$$ 
one has the following: 
\begin{equation}\label{Deltaak}
\Delta A_{k+1}=k U_{k+1}+U_{k+2}(\Delta\otimes Id^{\otimes k}).
\end{equation}
\rm

\medskip


\noindent{\sl Proof--} The proof is by induction on $k$. The relation
(\ref{dlaw}) reads:
$$\Delta\mu=Id\otimes Id+(Id\otimes\mu)(\Delta\otimes Id)+(\mu\otimes
Id)\tau_{23}(\Delta\otimes Id).$$
Since $A_1=Id,\ A_2=\mu$ and the domain is $Z_2$, one has:
$$\Delta A_2=U_2+(A_1\otimes A_2+A_2\otimes A_1)(\Delta\otimes Id)\ {\rm on}\  Z_2,$$
hence the relation (\ref{Deltaak}) is satisfied for $k=1$. For $k>0$
one has

\begin{align*}
     \Delta A_{k+1}= &\Delta\mu U_{k+1}    \\
      =&U_{k+1}+(Id\otimes\mu+(\mu\otimes Id)\tau_{23})(\Delta\otimes Id)U_{k+1}.
       \end{align*}

By induction one has
\begin{align*}
(\Delta\otimes Id)U_{k+1}=& \sum_{l=1}^k {k-1\choose l-1} \Delta A_l\otimes A_{k+1-l}\\
=& \sum_{l=2}^k {k-1\choose l-1} (l-1) U_{l}\otimes A_{k+1-l}\\
&+\sum_{l=1}^k {k-1\choose l-1} U_{l+1}(\Delta\otimes Id^{\otimes(l-1)})\otimes A_{k+1-l}
\end{align*}

$(Id\otimes\mu+(\mu\otimes Id)\tau_{23})$ applied to the first summand gives $(k-1)U_{k+1}$. 
Indeed:
\begin{multline*}
(Id\otimes \mu)\sum_{l=2}^k {k-1\choose l-1} (l-1) U_{l}\otimes A_{k+1-l}= \\
(Id\otimes\mu)\sum_{l=2}^k {k-1\choose l-1} (l-1)\sum_{j=1}^{l-1}
  {l-2\choose j-1}A_j\otimes A_{l-j} \otimes A_{k+1-l}=\\
\sum_{j=1}^{k-1} \sum_{1\leq l-j \leq k-j} (k-1){k-2\choose
  j-1}{k-j-1\choose l-j-1} A_j\otimes \mu( A_{l-j} \otimes A_{k+1-l})=\\
(k-1)\sum_{j=1}^{k-1} {k-2\choose j-1}A_j\otimes A_{k+1-j},
\end{multline*}
and
\begin{multline*}
(\mu\otimes Id)\tau_{23}\sum_{l=2}^k{k-1\choose l-1} (l-1) U_{l}\otimes A_{k+1-l}= \\
(\mu\otimes Id)\tau_{23}\sum_{l=2}^k {k-1\choose l-1} (l-1)\sum_{j=1}^{l-1}{l-2\choose j-1}A_j\otimes A_{l-j} \otimes A_{k+1-l}  \  \stackrel{{\rm on}\ Z_{k+1}}{=} \\
(\mu\otimes Id)\sum_{l=2}^k {k-1\choose l-1} (l-1)\sum_{j=1}^{l-1}{l-2\choose j-1}A_j\otimes A_{k+1-l} \otimes A_{l-j} \ \  \stackrel{\scriptstyle (l=k+1+j-\alpha)}{=} \\
\sum_{\alpha=2}^k\sum_{j=1}^{\alpha-1}{k-1\choose
  k+j-\alpha}(k+j-\alpha){k+j-\alpha-1\choose j-1} \mu(A_j\otimes A_{\alpha-j})\otimes A_{k+1-\alpha}= \\
\sum_{\alpha=2}^k(k-1) {k-2\choose \alpha-2} A_\alpha\otimes A_{k+1-\alpha}.
\end{multline*}

Adding the two results, one gets $(k-1)U_{k+1}$, thus
\begin{multline*}
\Delta A_{k+1}=kU_{k+1}+(Id\otimes\mu+(\mu\otimes Id)\tau_{23})\sum_{l=1}^k {k-1\choose l-1} U_{l+1}(\Delta\otimes Id^{\otimes(l-1)})\otimes A_{k+1-l}.
\end{multline*}

The computation of the second summand is similar to the previous one, thus on $Z_{k+1}$ one has
\begin{multline*}
(Id\otimes\mu+(\mu\otimes Id)\tau_{23})\sum_{l=1}^k {k-1\choose l-1} U_{l+1}(\Delta\otimes Id^{\otimes(l-1)})\otimes A_{k+1-l}=\\
U_{k+2}(\Delta\otimes Id^{\otimes k}).\  \Box
\end{multline*}


\subsection{Proof of property (1) of the fundamental lemma} For $x\in H$, since 
$\Delta^k(x)\in Z_{k+1}$ the formula (\ref{Deltaak}) can be applied:

\begin{multline*}
\Delta(e(x))=\Delta(x)+ \sum_{k\geq 1} \frac{(-1)^k}{k!}\Delta A_{k+1}\Delta^k(x) \\
\stackrel{\rm{by\ } (\ref{Deltaak})}{=}   \Delta(x)+ \sum_{k\geq 1} \frac{(-1)^k}{k!} kU_{k+1}\Delta^k(x)+ \sum_{k\geq 1}
 \frac{(-1)^k}{k!} U_{k+2}\Delta^{k+1}(x) \\
=\Delta(x)-U_2\Delta(x)=0. 
\end{multline*}
Hence $\Ima(e)\subset \Prim(H)$. Furthermore, if $y\in \Prim(H)$ then $e(y)=y$, so $e$ is an 
idempotent.  \hfill $\Box$


\subsection{Lemma}\it Let $H$ be a pre-Lie algebra and a \nap\  coalgebra satisfying the relation
(\ref{dlaw}). We denote by $\circ$ the product on $H$ as well as the action of $H$ on $H^{\otimes k}$
defined in \ref{defaction}. For $y\in H$, we define a map $\delta_y:H^{\otimes k}\rightarrow H^{\otimes(k+1)}$
by 
\begin{equation}\label{my}
\delta_y(x_1\otimes\ldots\otimes x_k)=\sum_{i=1}^k x_1\otimes\ldots x_i\otimes y\otimes x_{i+1}\otimes\ldots
\otimes x_k.
\end{equation}
The following equation holds
\begin{equation}\label{deltaitere}
\Delta^k(x\circ y)=\Delta^{k}(x) \circ y+\delta_y(\Delta^{k-1}(x)),\quad \forall k\geq 1
\end{equation}
\rm

\medskip


\noindent{\sl Proof--} The proof is by induction on $k$. For $k=1$ it is the relation (\ref{dlaw}).
For $k>1$ one has

\begin{multline*}
\Delta^{k+1}(x\circ y)=(\Delta\otimes Id^{\otimes k})(\Delta^k(x\circ y))
=(\Delta\otimes Id^{\otimes k})(\Delta^k(x)\circ y+\delta_y(\Delta^{k-1}(x))\\
=\Delta^{k+1}(x)\circ y-x_{(1)}\circ y\otimes x_{(2)}\otimes\ldots\otimes x_{(k+2)} 
-x_{(1)}\otimes x_{(2)}\circ y\otimes\ldots\otimes x_{(k+2)}\\
+\Delta(x_{(1)}\circ y)\otimes x_{(2)}\otimes\ldots
\otimes  x_{(k+1)} 
+\delta_y(\Delta^k(x))-x_{(1)}\otimes y\otimes x_{(2)}\otimes\ldots x_{(k+1)}\\
= \Delta^{k+1}(x)\circ y+\delta_y(\Delta^k(x)).\hspace{2cm} \Box
\end{multline*}


\subsection{Proposition}\it Let $H$ be a pre-Lie algebra and a \nap\  coalgebra satisfying
relation (\ref{dlaw}). The following equality holds, for $k\geq 0$, 
$x\in (H^{\otimes(k+1)})^{\S_1\times \S_k}$ and  $y\in H$:
\begin{equation}
\label{eqderak}
(k+1)A_{k+1}(x\circ y)=A_{k+2}\delta_y(x).
\end{equation}\rm

\medskip

\noindent{\sl Proof--} The proof is by induction on $k$. For $k=0$, the equality reads
$A_1(x\circ y)=x\circ y=\mu\delta_y(x)=A_2\delta_y(x).$  Denote by $m_y:H^{\otimes k}
\rightarrow H^{\otimes k}$ the map which associates $x\circ y$ to
$x\in H^{\otimes k}$. We
denote sometimes $m_y$ and $\delta_y$ by $m_y^k$ and $\delta_y^k$ if
$k$ is not fixed. For instance we have the formulas, 
for $\bar x=x_1\otimes\ldots\otimes x_{k+1}\in H^{\otimes (k+1)},\ y\in H$ and $1\leq l\leq k+1$:

\begin{equation}\label{relm-delta}
\begin{split}
\delta_y^{k+1}(\bar x)=& (\delta_y^{l-1}\otimes 
Id^{\otimes(k+2-l)})(\bar x)+x_1\otimes\ldots\otimes x_l\otimes y\otimes x_{l+1}\otimes\ldots\otimes x_{k+1}\\
&+ (Id^{\otimes l}\otimes \delta_y^{k+1-l})(\bar x)\\
m_y^{k+1}(\bar x)=&(m_y^{l}\otimes 
Id^{\otimes(k+1-l)})(\bar x)+(Id^{\otimes l}\otimes m_y^{k+1-l})(\bar x)
\end{split}
\end{equation}

Hence the formula (\ref{eqderak}) reduces to
$(k+1)A_{k+1}m^{k+1}_y=A_{k+2}\delta^{k+1}_y$
on the invariant space 
$(H^{\otimes (k+1)})^{\S_1\times \S_k}$.

\begin{align*}
A_{k+2}\delta_y^{k+1}(\bar x)=&
\mu (\sum_{l=1}^{k+1} {k\choose l-1} A_l\otimes A_{k+2-l})\delta_y^{k+1}(\bar x) \\
\stackrel{{\rm by\ }(\ref{relm-delta})}{=}&
\mu( \sum_{l=2}^{k+1} {k\choose l-1} A_l\delta_y^{l-1}\otimes A_{k+2-l}+
\sum_{l=1}^{k} {k\choose l-1} A_l\otimes A_{k+2-l}\delta_y^{k+1-l})(\bar x)  \\
&+\mu \sum_{l=1}^{k+1} {k\choose l-1} (A_l\otimes A_{k+2-l})(x_1\ldots x_l\otimes y\otimes x_{l+1}\ldots x_{k+1})\\
\stackrel{\rm{by\ induction}}{=}&
\mu[\sum_{l=2}^{k+1} {k\choose l-1}(l-1) A_{l-1}m_y^{l-1}\otimes A_{k+2-l}\\
&+\sum_{l=1}^{k} {k\choose l-1} (k+1-l)A_l\otimes A_{k+1-l}m_y^{k+1-l}](\bar x)  \\
&+\mu \sum_{l=1}^{k+1} {k\choose l-1} (A_l\otimes A_{k+2-l})(x_1\ldots x_l\otimes y\otimes x_{l+1}\ldots x_{k+1})\\
\end{align*}

But the first two terms of the right hand side, after a change of variables and thanks to the equality
$l{k\choose l}=(k+1-l){k\choose l-1}=k {k-1\choose l-1}$ writes

$$\mu\sum_{l=1}^{k} k{k-1\choose l-1}(A_{l}m_y^{l}\otimes A_{k+1-l}+
A_l\otimes A_{k+1-l}m_y^{k+1-l})(\bar x) 
=kA_{k+1}m_y(\bar x).$$
\bigskip

Again, by induction on $k$ we prove the following formula

\begin{equation}
\label{petitdernier}
\sum_{l=1}^{k+1} {k\choose l-1} A_l(x_1\ldots x_l)\circ  A_{k+2-l}(y\otimes x_{l+1}\ldots x_{k+1})=A_{k+1}m_y(\bar x).
\end{equation}

If $k=0$, this equation reads $x_1\circ y=m_y(x_1).$ If $k=1$ it reads
$x_1\circ(y\circ x_2)+(x_1\circ x_2)\circ y=x_1\circ (x_2\circ y)+(x_1\circ y)\circ x_2$ which is exactly the
pre-Lie relation. For $k>1$ one has

\begin{align*}
A_{k+1}m_y(\bar x)=&\mu(\sum_{l=1}^k{k-1\choose l-1} (A_lm_y\otimes A_{k+1-l}+
A_l\otimes A_{k+1-l}m_y))(\bar x) \\
\stackrel{{\rm by\ induction}}{=}\ &\sum_{l=1}^k\sum_{j=1}^l {k-1\choose l-1}{l-1\choose j-1}
(A_j(x_1\ldots x_j)\circ A_{l+1-j}(y\otimes x_{j+1}\ldots x_l))\circ \\
&A_{k+1-l}(x_{l+1}\ldots x_{k+1}) \\
&+\sum_{j=1}^k\sum_{u=1}^{k+1-j} {k-1\choose j-1}{k-j\choose u-1}
\underbrace{A_j(x_1\ldots x_j)}_{R}\circ \\
&(\underbrace{A_{u}(x_{j+1}\ldots x_{j+u})}_{T}\circ 
\underbrace{A_{k+2-j-u}(y\otimes x_{j+u+1}\ldots x_{k+1})}_{S}) 
\end{align*}

A change of variables in the first summand ($(j,l)\mapsto (j,u)$
with $l+1-j=k+2-j-u$),  gives a coefficient 
${k-1\choose k-u}{k-u\choose j-1}={k-1\choose j-1}{k-j\choose u-1}$, 
which is the same as the coefficient in the second summand. Furthermore since $\bar x\in (H^{\otimes (k+1)})^{\S_1\times \S_k}$
one has
\begin{multline*}
x_1\ldots x_j\otimes x_{j+1}\ldots x_{j+(k+1-j-u)}\otimes
\ldots\otimes x_{k+1}=\\
x_1\ldots x_j\otimes  x_{j+u+1}\ldots x_{k+1}\otimes x_{j+1}\ldots
x_{j+u}.
\end{multline*} 
The right hand side of the equality has the form $(R\circ S)\circ T+R\circ (T\circ
S)$, hence by the pre-Lie relation it writes also  $R\circ (S\circ T)+(R\circ T)\circ S+$ which gives the following:
 
 \begin{align*}
A_{k+1}m_y(\bar x)=&\sum_{j=1}^k\sum_{u=1}^{k+1-j}
	 {k-1\choose j-1}{k-j\choose u-1}[
A_j(x_1\ldots x_j)\circ \\
&(A_{k+2-j-u}(y\otimes x_{j+u+1}\ldots x_{k+1})\circ A_{u}(x_{j+1}\ldots x_{j+u})) \\
&+(A_j(x_1\ldots x_j)\circ A_{u}(x_{j+1}\ldots x_{j+u}))\\ 
&\circ A_{k+2-j-u}(y\otimes x_{j+u+1}\ldots x_{k+1})] 
\end{align*}
\begin{align*}
\stackrel{{\rm invariance\ of\ } \bar
  x}{=}&\sum_{j=1}^k\sum_{u=1}^{k+1-j} {k-1\choose j-1}{k-j\choose u-1}[
A_j(x_1\ldots x_j)\circ \\
&(A_{k+2-j-u}(y\otimes x_{j+1}\ldots x_{k+1-u})\circ A_{u}(x_{k+2-u}\ldots x_{k+1})) \\
&+(A_j(x_1\ldots x_j)\circ A_{u}(x_{j+1}\ldots x_{j+u}))\\
&\circ A_{k+2-j-u}(y\otimes x_{j+u+1}\ldots x_{k+1})] \\
\end{align*}

Since ${k-1\choose j-1}{k-j\choose u-1}={k-j\choose k+1-j-u}{k-1\choose j-1}$ the first summand gives
$$\sum_{j=1}^{k}{k-1\choose j-1} A_j(x_1\ldots x_j)\circ A_{k+2-j}
(y\otimes x_{j+1}\ldots x_{k+1}).$$ 
Since ${k-1\choose j-1}{k-j\choose u-1}={j+u-2\choose j-1}{k-1\choose j+u-2}$ the second summand gives
$$\sum_{\alpha=2}^{k+1}{k-1\choose \alpha-2} A_\alpha(x_1\ldots x_{\alpha})\circ A_{k+2-\alpha}
(y\otimes x_{\alpha+1}\ldots x_{k+1}).$$ 
Adding the two results one gets the relation (\ref{petitdernier}).
\hfill$\Box$


\subsection{Proof of the property (2) of the fundamental lemma} Let $x,y$ be in $H$. Since
$\Delta^k(x)\in (H^{\otimes(k+1)})^{\S_1\times \S_k}$ one can apply formula (\ref{eqderak}).

\begin{align*}
e(x\circ y)=& x\circ y+ \sum_{k>0} \frac{(-1)^k}{k!} A_{k+1}\Delta^k(x\circ y) \\
\stackrel{(\ref{deltaitere})}{=}& x\circ y+\sum_{k>0}
\frac{(-1)^k}{k!} A_{k+1}(\delta_y(\Delta^{k-1}(x))+\Delta^k(x)\circ y)\\
\stackrel{(\ref{eqderak})}{=}& x\circ y+\sum_{k>0}
\frac{(-1)^k}{k!} k A_{k}(\Delta^{k-1}(x)\circ y)+\sum_{k>0}
\frac{(-1)^k}{k!} A_{k+1}(\Delta^{k}(x)\circ y)\\
=&x\circ y-A_1(\Delta^0(x)\circ y)=0 \hspace{6cm} \Box
\end{align*}


 \section{Using cogroups}
 
 In \cite{Frbull} and \cite{Fr98} B. Fresse proves that a cogroup in the category of connected graded
 $\mathcal P$-algebras, where $\mathcal P$ is an operad is a free
 $\mathcal P$-algebra. In \cite{O99} J.-M. Oudom
 generalises this theorem to comagma in this category. This is a
 generalisation of a theorem of Leray \cite{Le45} for commutative
 algebras and of Berstein \cite{B65} 
for associative algebras. This technics was also studied by
R. Holtkamp, for instance in the context of $\mathcal P$-Hopf operads
(see e.g. his habi\-li\-tationsschrift \cite{Holt04}). The purpose of this section is to present
 the theorems of rigidity for associative algebras and for pre-Lie
 algebras by the way of abelian cogroups in some categories.
 
 \medskip
 
 \subsection{Coproduct for algebras over an operad}\label{coproduct}
 We recall the notation of Fresse. Let ${\mathcal O}$ be an operad. The category
 of graded algebras over ${\mathcal O}$ is endowed with a {\sl coproduct} $\vee$: let $A$ and $B$ be
 two ${\mathcal O}$-algebras then $A\vee B={\mathcal O}(A\oplus B)/rel$ where $rel$ is
 generated by relations of the form 
$$\mu\otimes (x_1\otimes\ldots\otimes x_n)-\mu(x_1,\ldots,x_n)$$ 
for all  $\mu\in{\mathcal O}(n)$ and $x_1\otimes\ldots\otimes x_n \in A^{\otimes n}\cup
 B^{\otimes n}$. In particular $A$ and $B$ are a summand of $A\vee
 B$. There is a natural morphism $T:A\vee B\rightarrow B\vee A$. 
A cogroup in this category is a connected graded ${\mathcal O}$-algebra together with a map
 $\phi: A\rightarrow A\vee A$ such that $\phi$ is coassociative and is
 the identity on each copy of $A$. A cogroup is abelian or
 cocommutative if $T\phi=\phi$. The theorem of Fresse we use is the following

\medskip

{\bf Theorem}\cite{Fr98}. \it  Fix a ground field of characteristic 0. let
$\mathcal P$ be a unital operad. If $R$ is a connected graded
$\mathcal P$-algebra equipped with a cogroup structure, then $R$ is a free
graded $\mathcal P$-algebra.\rm

 \subsection{Theorems of Berstein \cite{B65}  and Loday-Ronco  \cite{LR04}
   revisited} Since we
 work in the context of operads, we are concerned with non-unital
 associative algebras. In this context, a {\sl non-unital infinitesimal bialgebra} is
 a non-unital associative algebra together with a noncounital coassociative coproduct $\Delta$ satisfying the relation
 \begin{equation}\label{uilaw}
 \Delta(x\cdot y)=x\otimes y+ x_{(1)}\otimes x_{(2)}\cdot y+x\cdot y_{(1)}\otimes y_{(2)}
 \end{equation}

 \noindent{\bf Theorem--}\it Let $H$ be an object in the category of
 graded connected non-unital
associative algebras. The following are equivalent
 \begin{itemize}
 \item[a)] $H$ is a non-unital infinitesimal bialgebra,
 \item[b)] $H$ is free,
  \item[c)] $H$ is an abelian cogroup. 
\end{itemize}
 \rm
 
 \medskip
 
 The equivalence a) and b) is a theorem by Loday
 and  Ronco  \cite{LR04} on the structure of unital infinitesimal bialgebras 
 adapted to the non-unital case. The equivalence b) and c) is a consequence
 of a theorem by Berstein \cite{B65}. Indeed the equivalence a) and c) can be proved independently.
 Here is the proof: 

\smallskip

Recall that the coproduct in the category of associative algebras is given by
 $H_1\star H_2=H_1\oplus H_2\oplus \bigoplus_n \underbrace{H_1\otimes
 H_2\otimes\ldots}_{alternative\ tensor\ n\ times}
 \oplus \underbrace{H_2\otimes H_1\otimes\ldots}_{n\ times}$, hence
 $$H\star H=\bigoplus_{n\geq 1} (1,H^{\otimes n})+(2,H^{\otimes n}),$$ 
 where the label 1 or 2 indicates the beginning of the alternative tensor.
 The product of two elements
 $(i,a^n)$ and $(j,b^m)$ with $a^n\in H^{\otimes n}$ and $b^m\in H^{\otimes m}$ lies in the copy $i$ of $H$ and is
 the concatenation of $a$ and $b$ if the last term of $a^n$ lies in the copy $\not=j$ of $H$
 and is $a^n_1\otimes\ldots\otimes a^n_{n-1}\otimes a^n_n\cdot b^m_1\otimes\ldots\otimes b^m_m$
 otherwise. If $H$ is an abelian cogroup in the category of associative algebras then there is
 a map $\phi: H\rightarrow H\star H$ which is coassociative and is the identity on H:
 $\phi(h)=(1,h)+(2,h)+(1,h_1\otimes h_2)+(2,h_1\otimes h_2)+\cdots $. Define $\Delta:H\rightarrow H\otimes H$
 by $\Delta(h)=h_1\otimes h_2$ the projection $\pi$ of $\phi$ onto the
 vector space $(1,H^{\otimes 2})$. The coassociativy of $\phi$ implies the one of $\Delta$. Since $\phi$ is a morphism
of associative algebras, one has
\begin{multline*}
\Delta(a\cdot b)=\pi\phi(a\cdot b)=\\
\pi[((1,a)+(1,a_1\otimes a_2)+\ldots )\cdot((1,b)+(2,b)+(1,b_1\otimes b_2)+\ldots )]\\
=a\otimes b+a\cdot b_1\otimes b_2+a_1\otimes a_2\cdot b,
\end{multline*}
which is the relation (\ref{uilaw}).
\smallskip

To prove the converse it is easy to check that the map
$\phi:H\rightarrow H\star H$ defined by $\phi(a)=\sum_{n\geq 1} (1,\Delta^{n-1}(a))+
(2,\Delta^{n-1}(a))$ is well defined (the algebra is graded
connected), is coassociative, cocommutative and
is a morphism of algebras thanks to the relation

\begin{multline*}
\Delta^k(a\cdot b)=\sum_{i=1}^k a_{(1)}\otimes\ldots\otimes a_{(i)}\otimes b_{(1)}\otimes\ldots\otimes
b_{(k+1-i)} \\
+\sum_{i=1}^{k+1} a_{(1)}\otimes\ldots\otimes a_{(i)}\cdot b_{(1)}\otimes\ldots\otimes
b_{(k+2-i)}. \hspace{0.5cm} \Box
\end{multline*}

 \subsection{Cogroups in the category of pre-Lie algebras} Let $T$ be a rooted tree. A {\sl leaf} of $T$ is
 a vertex with no incoming edges. For $v\in\Vertex(T)$ and for some
 rooted trees $S_1,\ldots S_n$ the tree $T\circ_v\{S_1,\ldots, S_n\}$
 is the tree obtained by grafting the trees $S_i$ on the vertex $v$ of $T$.

\medskip

Denote by $\PL$ the operad defining pre-Lie algebras. Recall that
 $\PL(V)$ denotes the free pre-Lie algebra generated by $V$. 
Let $L_1$ and $L_2$ be two pre-Lie algebras. The coproduct of $L_1$ and $L_2$
 is $L_1\vee L_2=\PL(L_1\otimes L_2)/rel$ where the equivalence relation is the one defined in \ref{coproduct}.
 It is not difficult to prove the following lemma
 
 \subsubsection{Lemma}\it Let $L$ be a pre-Lie algebra, and $(e_\alpha)_{\alpha}$ be a basis
 of its  underlying vector space. Then a basis of $L\vee L$ is given by rooted trees labeled by
 the set $(i,e_\alpha)_{i\in\{1,2\},\alpha}$ satisfying the ``leaf condition'': the labeling $i$ of any
 leaf is different from the one of its adjacent vertex.\rm
 
 \medskip
 
 Furthermore, its pre-Lie structure is given by the following: let
 $T,S$ be two trees. If $deg(S)>1$ then
 $T\circ S=\sum_{v\in\Vertex(T)} T\circ_v\{S\}$. 
If $deg(S)=1$, assume $S$ is a single vertex
 labeled by $(1,e_\alpha)$ for instance. If $v\in\Vertex(T)$ is
 labeled by 2, set $T\hat\circ_v S=T\circ_v \{S\}$. If $v$ is
 labeled by $1$, then $T\circ_v \{S\}$ does not
 respect the ``leaf condition''. Denote by $T_1^v,\ldots T_k^v$ the incoming subtrees
 of $T$ at the vertex $v$; denote by $T_{v\leftarrow w}$ the tree obtained by replacing $v$
 by $w$. If $v=(1,e_\beta)$, we define
 $$T\hat\circ_v (1,e_\alpha)=T_{v\leftarrow (1,e_\beta\circ e_\alpha)}-
 \sum_{I,J} T\circ_v \{T_I,(1,e_\alpha)\circ T_J\}.$$
 where the sum is over all partition $I\cup J$ of $\{1,\ldots,k\}$
 such that $J\not=\emptyset.$ Then for $deg(S)=1$,
$T\circ S=\sum_{v\in\Vertex(T)} T\hat\circ_v S.$

 \subsubsection{Lemma}\it If $L$ is a pre-Lie algebra together with a cogroup structure 
 $\phi: L\rightarrow L\vee L$ then $L$ is endowed with a \nap\  coproduct $\Delta$ satisfying
 the relation (\ref{dlaw}). \rm
 
 \medskip
 
 \noindent{\sl Proof--} The relation $(Id-\tau_{23})(\Delta\otimes
 Id)\Delta$ is a consequence of the associativity of $\phi$. 
The relation (\ref{dlaw}) is a consequence of $\phi$ being a morphism of pre-Lie algebras.\hfill$\Box$
 
 \subsubsection{Theorem}\it Let $L$ be an object in the category of graded connected pre-Lie  algebras.
 The following are equivalent
 \begin{itemize}
 \item[a)] $L$ is a \nap\  coalgebra satisfying relation (\ref{dlaw})
 \item[b)] $L$ is free. 
  \item[c)] $L$ is an abelian cogroup. 
\end{itemize}
 \rm

\medskip

\noindent{\sl Proof--} The equivalence a) and b) is the purpose of our paper. The equivalence b) and c)
is the theorem of Fresse applied to the operad preLie. c) implies a) from the previous lemma. \hfill $\Box$
 
 \vspace{1cm}

\noindent{\bf Aknowledgments.} I would like to thank J.-L. Loday for
fruitful discussions and J.-M. Oudom for giving me details on
abelian cogroups. I'm indebted to the referee for his comments and
suggestions.

\bibliographystyle{amsplain} 
\bibliography{biblast.bib}

\bigskip

\end{document}